
\magnification\magstep1
\hsize 30pc
\vsize 44pc

\font\bfab = cmbx9
\font\fa=cmr17
\font\fb=cmr12

\font\fab=cmr9
\font\sfab=cmr7
\font\fs=cmr6
\font\fd=cmr5
\font\slab=cmsl9
\font\tab= cmmi9
\font\sab= cmmi7
\font\ssab= cmmi6
\font\abst= cmsy9
\font\sabs= cmsy7
\font\ssabs= cmsy6
\font\itab= cmti9

\def\figfont{
    \textfont0 = \fab
    \scriptfont0 = \sfab
    \scriptscriptfont0 = \fs
    \textfont1 = \tab
    \scriptfont1 = \sab
    \scriptscriptfont1 = \ssab
    \textfont2 = \abst
    \scriptfont2 = \sabs
    \scriptscriptfont2 = \ssabs
    \let \sl = \slab
    \let \bf = \bfab
    \let \it = \itab
    \baselineskip 9pt
    \fab}

\def\leftheadline{\ifnum\pageno=\count100 \hfill%
  \else\rm\folio\hfil\it\shortauthor\fi}
\def\rightheadline{\ifnum\pageno=\count100 \hfill%
  \else\it\shorttitle\hfil\rm\folio\fi}


\def\CC{{\rm C\kern-.18cm\vrule width.6pt height 6pt depth-.2pt
\kern.18cm}}

\def\NN{{\mathop{{\rm I}\kern-.2em{\rm N}}\nolimits}}

\def\PP{{\mathop{{\rm I}\kern-.2em{\rm P}}\nolimits}}

\def\RR{{\mathop{{\rm I}\kern-.2em{\rm R}}\nolimits}}

\def\RRt{{\fa I}\kern-.2em{\fa R}}

\def\ZZ{{\mathop{{\rm Z}\kern-.28em{\rm Z}}\nolimits}}





\def\makebold#1{\mathord{\setbox0=\hbox{#1}%
       \copy0\kern-\wd0%
       \raise\dimen1\copy0\kern-\wd0%
       {\advance\dimen1 by \dimen1\raise\dimen1\copy0}\kern-\wd0%
       \kern\dimen0\raise\dimen1\copy0\kern-\wd0%
       {\advance\dimen1 by \dimen1\raise\dimen1\copy0}\kern-\wd0%
       \kern\dimen0\raise\dimen1\copy0\kern-\wd0%
       {\advance\dimen1 by \dimen1\raise\dimen1\copy0}\kern-\wd0%
       \kern\dimen0\raise\dimen1\copy0\kern-\wd0%
       \kern\dimen0\box0}}


\def\cro{\cr \noalign{\vskip 4pt}}


%

\def\frac#1#2{{#1 \over #2}}


\def\noin{\noindent}










\def\pf{\noindent{\bf Proof: }}

\def\eop{\makeblanksquare6{.4}}

\def\makeblanksquare#1#2{
\dimen0=#1pt\advance\dimen0 by -#2pt
      \vrule height#1pt width#2pt depth0pt\kern-#2pt
      \vrule height#1pt width#1pt depth-\dimen0 \kern-#1pt
      \vrule height#2pt width#1pt depth0pt \kern-#2pt
      \vrule height#1pt width#2pt depth0pt
}



\def\abstract#1{\bigskip\bigskip\medskip%
    {\narrower \baselineskip 9pt \fab \noindent {\bf Abstract.~~}%
    \textfont0 = \fab
    \scriptfont0 = \sfab
    \scriptscriptfont0 = \fs
    \textfont1 = \tab
    \scriptfont1 = \sab
    \scriptscriptfont1 = \ssab
    \textfont2 = \abst
    \scriptfont2 = \sabs
    \scriptscriptfont2 = \ssabs
    \let \it = \itab
    \let \sl = \slab
    #1\bigskip}\medskip}
\def\author#1{\bigskip\bigskip\centerline{\fb #1}}

\def\copyright{\hbox{{\fb o}\kern-.61em \raise .46ex \hbox{\fd c}}}
\def\title#1{\centerline {\fa #1}}
\def\titwo#1{\medskip \centerline {\fa #1}}

\def\titexp#1#2{\hbox{{\fa #1} \kern-.25em \raise .90ex \hbox{\fb #2}}\/}
\def\titsub#1#2{\hbox{{\fa #1} \kern-.25em \lower .60ex \hbox{\fb #2}}\/}


\def\sect#1#2{\goodbreak\bigskip\centerline{\bf\S #1. #2}\medskip
    \noindent\ignorespaces\secnum=#1\eqnum=0\proclaimnum=0}





\def\Remark#1.{\goodbreak\medskip\noin {\bf Remark#1.}}



\def\ref{\smallskip\global\advance\refnum by 1 \item{\the\refnum.}}
\newcount\refnum \refnum = 0




\newcount\blackmarks\blackmarks=0
\newcount\eqnum\eqnum = 0
\newcount\figurenum\figurenum=0
\newcount\proclaimnum\proclaimnum=0
\newcount\refnum\refnum = 0
\newcount\tablenum\tablenum=0
\newcount\secnum\secnum = 0

 \newread\testfl
 \def\inputifthere#1{\immediate\openin\testfl=#1
    \ifeof\testfl\message{(#1 does not yet exist)}
    \else\input#1\fi\closein\testfl}

 \inputifthere{\jobname.aux}
 \newwrite\mem
 \immediate\openout\mem=\jobname.aux

\def\plazieres{\expandafter\ifx\csname\griff\endcsname\relax%
  \xdef\esfehlt{\griff}\blackmark\else{\csname\griff\endcsname}\fi}

\def\definieres{\expandafter\xdef\csname\griff\endcsname{\inhalt}%
 \def\blankkk{ }\expandafter\immediate\write\mem{%
 \string\expandafter\def\string\csname%
 \blankkk\griff\string\endcsname{\inhalt}}}

\def\blackmark{\ifnum\blackmarks=0\global\blackmarks=1%
 \write16{===========================================================}%
 \write16{Some forward reference not yet defined. Re-TeX this file!}%
 \write16{===========================================================}%
 \fi\immediate\write16{undefined forward reference: \esfehlt}%
 {\vrule height10pt width2pt depth2pt}\esfehlt%
 {\vrule height10pt width2pt depth2pt}}

\def\showlabel#1{\marginal{#1}}
\def\marginal#1{\strut\vadjust{\kern-\strutdepth%
\vtop to \strutdepth{\baselineskip\strutdepth\vss\llap{\fiverm#1\ }\null}}}
\def\strutdepth{\dp\strutbox}


\newif\ifdraft

\newcount\hour\newcount\minutes
\def\draft{\drafttrue
\def\comment##1{{\bf comment: ##1}}
\headline={\sevenrm \hfill
\ifx\filenamed\undefined\jobname\else\filenamed\fi%
(.tex) (as of \ifx\updated\undefined???\else\updated\fi)
 \TeX'ed at {\hour\time\divide\hour by 60{}%
\minutes\hour\multiply\minutes by 60{}%
\advance\time by -\minutes
\the\hour:\ifnum\time<10{}0\fi\the\time\  on \today\hfill}}
}

\def\today{\number\day\space%
\ifcase\month\or January\or February\or March\or April\or May\or June\or
 July\or August\or September\or October\or November\or December\fi%
\space\number\year}

\def\en#1{%
  \global\advance\eqnum by 1
  \edef\griff{en#1}\edef\inhalt{\the\secnum.\the\eqnum}%
  \definieres\eqno(\inhalt)\ifdraft\rlap{\fiverm #1}\fi}

\def\enn#1{\global\advance\eqnum by 1
  \edef\griff{en#1}\edef\inhalt{\the\secnum.\the\eqnum}%
    \definieres(\inhalt)\ifdraft\rlap{\fiverm #1}\fi}


\def\fn#1{\global\advance\figurenum by 1
  \edef\griff{fn#1}\edef\inhalt{\the\figurenum}\definieres\inhalt
  \ifdraft\showlabel{#1}\fi}%


\def\pn#1{~\global\advance\proclaimnum by 1
  \edef\griff{pn#1}\edef\inhalt{\the\secnum.\the\proclaimnum}%
  \definieres{\inhalt}\ifdraft\showlabel{#1}\fi}


\def\rhl#1{\global\advance\refnum by 1
    \edef\griff{rn#1}\edef\inhalt{\the\refnum}\definieres%
    \wo{}\wo{\noexpand\smallskip\noexpand\item{\the\refnum.}}} 

\def\rn#1{\smallskip\global\advance\refnum by 1
   \edef\griff{rn#1}\edef\inhalt{\the\refnum}\definieres%
   \ifdraft\showlabel{#1}\vskip-\baselineskip\fi\item{\the\refnum.}}%

\def\cit#1{\edef\griff{rn#1}\plazieres}

\def\cite#1{[\cit{#1}]} 

\def\tn#1{\global\advance\tablenum by 1
   \edef\griff{tn#1}\edef\inhalt{\the\tablenum}\definieres \inhalt
   \ifdraft\showlabel{#1}\fi}%



\input epsf

%

%

\newtoks\lastname
\newtoks\firstname
\newtoks\au
\newtoks\aut
\newtoks\ti
\newtoks\tit
\newtoks\pb
\newtoks\pub
\newtoks\pl

\newif\ifonesofar
\def\concat#1{\edef\audef{{#1}}\au=\audef}
\def\decodeauthor#1, #2,#3;{\lastname={#1}\firstname={#2}%
\concat{\formfirstauthor}\onesofartrue%
\def\next{#3}\ifx\next\empty\else\decodemoreauthor#3;\fi}
\def\decodemoreauthor#1, #2,#3;{\lastname={#1}\firstname={#2}%
\def\next{#3}
\ifx\next\empty\let\formaut=\formlastauthor%
\ifonesofar\ifx\formotherauthor\undefined\else\let\formaut=\formotherauthor%
\fi\fi\concat{\the\au\formaut}%
\else\onesofarfalse\concat{\the\au\formnextauthor}\decodemoreauthor#3;\fi}

\def\refB #1; #2; #3 (#4); #5; {\decodeauthor#1,;%
   \ti={#2}\pb={#3}\pl={#4}\def\yr{#5}\formB}

\def\refD #1; #2; #3; #4; {\decodeauthor#1,;%
   \ti={#2}\pl={#3}\def\yr{#4}\formD}

\def\refJ #1; #2; #3; #4; #5; #6; {\decodeauthor#1,;%
    \ti={#2}\def\jr{#3}\def\vl{#4}\def\yr{#5}\def\pp{#6}\formJ}

\def\lookupp#1{{\global\aut={\vrule height15pt width15pt depth0pt}%
 \global\tit={{\bf the specified proceedings does not exist in our files}}%
 \xdef\edsop{}\global\pub={}\def#1{}\input \locbib proceed }}

\def\refproc #1(#2; #3; {\decodeproc#2; \xdef\yr{#3}}
\def\decodeproc#1), #2 (ed#3.),#4 (#5); {%
 \global\tit={#1}\global\aut={#2}\xdef\edsop{#3}\global
 \pub={#4}\global\pl={#5}}

\def\refP #1; #2; #3; #4; {\lookupp{#3}\decodeauthor#1,;%
        \ti={#2}\def\pp{#4}\formP}

\def\refQ #1; #2; (#3; #4; #5; {\decodeproc#3; \decodeauthor#1,;%
   \ti={#2}\def\yr{#4}\def\pp{#5}\formP}

\def\refR #1; #2; #3; #4; {\decodeauthor#1,;%
         \ti={#2}\def\is{#3}\def\yr{#4}\formR}

\def\refnew #1; {
 \wo{#1.}}


\newwrite\w
\def\wo#1{\immediate\write\w{#1}}
\def\startbib{\immediate\openout\w=\jobname.bib}
\def\endbib{\immediate\closeout\w \input\jobname.bib}
\newcount\refnum \refnum = 0

\def\locbib{}


%
\gdef\formfirstauthor{\the\lastname, \the\firstname}
\gdef\formnextauthor{, \the\firstname\the\lastname}
\gdef\formotherauthor{ and \the\firstname\the\lastname}
\gdef\formlastauthor{,\formotherauthor}
\gdef\formB{
\wo{\the\au,}
\wo{{\sl \the\ti,}}
\wo{\the\pb,}
\wo{\the\pl,}
\wo{\yr.}}
\gdef\formD{
\wo{\the\au,}
\wo{\the\ti,}
\wo{dissertation,}
\wo{\the\pl,}
\wo{\yr.}}
\gdef\formJ{
\wo{\the\au,}
\wo{\the\ti,}
\wo{\jr\noexpand\ {\noexpand\bf\vl}}
\wo{(\yr), \pp.}}
\gdef\formP{
\wo{\the\au,}
\wo{\the\ti,}
\wo{in}
\wo{{\sl \the\tit,}}
\wo{\the\aut\ (ed\edsop),}
\wo{\the\pub,}
\wo{\the\pl,}
\wo{\yr,}
\wo{\pp.}}
\gdef\formR{
\wo{\the\au,}
\wo{\the\ti,}
\ifx\is\empty\else\wo{\is,}\fi
\wo{\yr.}}




\overfullrule=0pt
\title{Factorization of}
\titwo{Multivariate Positive Laurent Polynomials}
\author{Jeffrey S. Geronimo\footnote{$^{1)}$}{\fab
School of Mathematics, The Georgia Institute of Technology, Atlanta, GA
30332, \hfill\break geronimo@math.gatech.edu. This author is partly supported by an NSF grant.} 
and Ming-Jun Lai\footnote{$^{2)}$}{\fab
Department of Mathematics, The University of Georgia, Athens, GA
30602, \hfill\break mjlai@math.uga.edu. This author is partly
supported by the National Science Foundation under grant EAR-0327577 and
the School of Mathematics, The Georgia Institute of Technology when he 
visited the School during the fall, 2004.}}

\abstract{Recently Dritschel proves that any positive multivariate Laurent polynomial can
be factorized into a sum of square magnitudes of polynomials. We first give another proof
of the Dritschel theorem. Our proof is based on the univariate matrix F\'ejer-Riesz theorem.
Then we discuss a computational method to find approximates of polynomial matrix factorization. 
Some numerical examples will be shown. Finally we discuss how to compute 
nonnegative Laurent polynomial factorizations in the multivariate setting.}

\sect{1}{Introduction}
We are interested in computing factorizations of nonnegative Laurent polynomials into sum of
squares of polynomials. That is, let 
$$
P(z)= \sum_{k=-n}^n p_k z^k
$$
be a Laurent polynomial, where $z=e^{i\theta}$. Suppose that $P(z)\ge 0$ for $|z|=1$. 
One would ask if there exists a polynomial $Q(z)=\displaystyle \sum_{k=0}^n q_kz^k$ such that
$$
P(z)= Q(z)^*Q(z), \eqno(1)
$$
where $Q(z)^*$ denotes the complex conjugate of $Q(z)$. This is the well-known Fej\'er-Riesz
factorization problem and it was resolved by Fej\'er [F'15] and by Riesz [R'15]. A natural
question is whether the results of Fej\'er and Riesz can be extended to the multivariate setting.
More generally, given a nonnegative multivariate trigonometric polynomial $P(z):=P(z_1,z_2,\cdots,
z_d)$ of coordinate degrees $\le n$, 
does there exist a finite number of polynomials $Q_k(z)$ such that  
$$
P(z)= \sum_k Q^*_k(z)Q_k(z), \eqno(2)
$$
i.e., can $P(z)$ be written as a sum of square magnitudes  (sosm) of polynomials. 
There is a vast amount of  literature related to the study of this problem 
and the results relevant to this paper may be summarized as follows:
\item{$1^\circ$} When $P(z)$ is nonnegative on 
the multi-torus $|z_1|=|z_2|=\cdots =|z_d|=1$ and the coordinate degrees of $Q_k$ are 
less than or equal to $n$, then the answer to the question is negative. (See 
[Calderon and Pepinsky'52] and [Rudin'63].)

\item{$2^\circ$} When $P(z)$ is strictly positive on the multi-torus and the coordinate degrees
of $Q_k$ are not specified, Dritschel has shown that 
the answer to the question is positive([Dritschel'04]). However the nonnegative case remains unresolved.  
\item{$3^\circ$} In the bivariate setting, Geronimo and Woerdeman gave a necessary and sufficient
condition in order for $P(z)=|Q(z)|^2$, where $Q(z)$ is a stable polynomial, i.e., $Q(z)\not=0$
inside and on the bi-torus ([Geronimo and Woerdeman'04]). 

\item{$4^\circ$} In the bivariate setting, there exist rational Laurent polynomials $Q_k(z)$ 
such that (2) holds.  Furthermore, $Q_k$ can be so chosen that the 
determinants of $Q_k$ containing only one variable Laurent  polynomials (cf. [Basu'01]).  

\item{$5^\circ$} In [McLean and Woerdeman'01], an algorithm was proposed to
find polynomials $P_k$ such that $P = \sum_k |P_k|^2$.  
The algorithm uses the so-called semi-definite programming.

Although the mathematical problem appears to be theoretical, it has many applications in engineering,
e.g., the design of autoregressive filters, construction of orthonormal wavelets (cf. [Daubechies'92]), 
construction of tight wavelet framelets (cf. [Lai and Stoeckler'04]), spectral estimation in
control theory (cf. [Sayed and Kailath'01]) and many other engineering applications mentioned in 
[McLean and Woerdeman'01].    
Thus, how to compute such factorization polynomials $Q_1, Q_2, \cdots,$ is  interesting and 
useful for applications. In this paper, we discuss a symobl approximation method 
studied in [Lai'94] for computing
such factorizations. The method was originally intended for
factorizing any nonnegative Laurent polynomials in the univariate setting. We use the ideas to 
give Dritschel's theorem another proof. The proof
provides a computational method to factor $P(z)$ into $Q_k(z)'s$. The paper
is organized as follows. In section 2, we first give  
a different proof of Dritschel's Theorem. A key in the proof is to factorize univariate Laurent
polynomial matrices. In section 3, we discuss how to compute the factorization
of positive Laurent polynomials matrices.  The method
used by [Lai'94] to compute approximate factorizations is extended to the
matrix case. Then in section 4, some numerical examples are computed
following the procedure in \S 2 and \S 3. Finally in section 5. the
nonegative case is considered.   

\sect{2}{Dritshel's Theorem}
We begin with reviewing the concept of the symbols of bi-infinite Toeplitz matrices and its properties.
For a given univariate Laurent polynomial $\displaystyle P(z)=\sum_{k=-n}^n p_k z^k$, 
we may view $P(z)$ as the symbol of a bi-infinite Toeplitz matrix  
${\cal P}:=(p_{i-j})_{i, j\in {\bf Z}}$. Indeed, for any absolutely summable
sequence ${\bf x}=(x_i)_{i\in {\bf Z}}$, i.e., $\displaystyle \sum_{i\in {\bf Z}}|x_i| <\infty$,
let $F({\bf x})= \sum_{j\in {\bf Z}}x_j z^j$ be the discrete Fourier transform (or z-transform)
of ${\bf x}$. Let ${\bf y}= {\cal P}{\bf x}$, then it is easy to see that
$$
F({\bf y})=  P(z) F({\bf x}).
$$
If the matrix ${\cal P}$ has a factorization ${\cal Q}$ which is a banded
upper triangular Toeplitz matrix such that 
$$
{\cal P}={\cal Q}^\dagger {\cal Q}, 
$$
the discrete Fourier transform of ${\bf y}={\cal Q} {\cal Q}^\dagger{\bf x}$ is 
$F({\bf y})=Q(z)^*Q(z)F({\bf x})$,
where ${\cal Q}^\dagger$ denotes the complex conjugate transpose of ${\cal Q}$. 
Thus, finding $P(z)= Q(z)^*Q(z)$ is equivalent to finding a banded upper triangular Toeplitz matrix ${\cal Q}$ such that
${\cal P}={\cal Q}^\dagger{\cal Q}$. 

It is easy to show that if $P(z)\ge 0$
for all $|z|=1$, then ${\cal P}$ is Hermitian and nonnegative definite. Clearly, ${\cal P}$ is
Hermitian since $P(z)$ is real. Furthermore for any absolutely summable sequence ${\bf x}$, 
we need to show
that ${\bf x}^\dagger  {\cal P}{\bf x}\ge 0$. Again writting ${\bf y}={\cal P}{\bf x}$, we know that
$${\bf x}^\dagger  {\bf y} = {1\over 2\pi}\int_0^{2\pi}\overline{F({\bf x})}F({\bf y})d\theta$$
where $z=e^{i\theta}$ and it follows that
$$
{\bf x}^\dagger  {\cal P}{\bf x}= \displaystyle {1\over 2\pi} \int_0^{2\pi} 
|F({\bf x})|^2 P(z)d\theta \ge 0 
$$
for any nonzero sequence ${\bf x}$. In particular, for 
$$
{\bf x}=(\cdots, 0,x_{-N},\cdots, x_0,
\cdots, x_N,0,\cdots)^T,
$$
the left-hand side in the above inequality gives ${\bf x}^\dagger  P_N {\bf x}$, where $P_N$ is a 
central section of ${\cal P}$. The above argument shows that $P_N$ is nonnegative definite. 

In the following
we will assume that $P(z)$ is strictly positive, in the sense that there exists a positive number
$\epsilon>0$ such that $P(z)\ge \epsilon$. When $P(z)$ is a matrix, we mean that 
$P(z)\ge \epsilon I$, where $I$ is the identity matrix of the same size as that of $P(z)$. 
When  $P(z)$ is strictly positive, we have
$$
{\bf x}^\dagger  {\cal P}{\bf x}= \displaystyle {1\over 2\pi} \int_0^{2\pi} |F({\bf x})|^2 
P(z)d\theta 
\ge \epsilon \|{\bf x}\|^2. 
$$ 
It follows that if $P(z)\ge \epsilon >0$, then $P_N\ge \epsilon>0$. 

We now consider  the factorization of multivariate Laurent 
polynomials.  Let us begin with a bivariate Laurent polynomial $P(z_1,z_2)$ first. 
That is, let 
$$
P(z_1,z_2)=\sum_{j=-n}^n \sum_{k=-n}^n p_{jk}z_1^j z_2^k\ge 0
$$
be a Laurent polynomial of coordinate degrees $\le n$. 
We would like to find a finite number of polynomials $Q_k$ such that 
$$
P(z_1,z_2)= \sum_k |Q_k(z_1,z_2)|^2.
$$
Denote by ${\bf z_1}=[1,z_1,z_1^2, \cdots, z_1^n]^T$ and write   
$$
P(z_1, z_2)= {\bf z_1}^\dagger\widetilde{P}(z_2) {\bf z_1}
$$
for a Hermitian matrix $\displaystyle \widetilde{P}(z_2) = \sum_{k=-n}^n
\tilde{p}_k z_2^k$, where each $p_k$ is an $(n+1)\times(n+1)$ Toeplitz
matrix. With a slight modification of an observation of  [McLean and Woerdeman'01, Theorem
2.1], we note that there are many ways to write $\widetilde{P}(z_2)$.  
If there is one $\widetilde{P}(z_2)$ 
which is  nonnegative definite then we can use the matrix F\'ejer-Riesz factorization 
(cf. e.g., in [Helson'64], [Mclean-Woerdeman'01], see also section 3) 
to 
find $\widetilde{Q}(z_2)$ such that 
$$
\widetilde{P}(z_2)= \widetilde{Q}^\dagger(z_2) \widetilde{Q}(z_2).
$$
That is,  we have
$$
P(z_1,z_2)=(\widetilde{Q}(z_2){\bf z_1})^\dagger\widetilde{Q}(z_2){\bf z_1}
$$  
which is clearly a sum of squares of polynomials.

The above discussion can be generalized to the multivariate 
setting and using an observation of [Dritschel'04] to the case that the
size of $\widetilde{P}(z_2)$ is larger than $(n+1)\times(n+1)$. 
For simplicity, let us consider  a trivariate Laurent polynomial $P(z_1,z_2,z_3)$ 
in $z_1=e^{i\theta_1}, z_2=e^{i\theta_2}, z_3=e^{i\theta_3}$ of coordinate degrees $\le n$.  
We first write $P(z_1,z_2,z_3)$ in a matrix format: 
$$
P(z_1,z_2,z_3)=\sum_{-n}^n p_k(z_2,z_3) z_1^k= {\bf z_1}^\dagger\widehat{P}(z_2,z_3) {\bf z_1}, 
$$
with
$$
{\bf z_1}=[1,z_1,\ldots, z_1^{m_1}]^T  \eqno(3)
$$
and $m_1\ge n$. There are many ways to write $\widehat{P}(z_2,z_3)$. To
capture this define the set of matrices
$$
{\cal F}(z_2,z_3)=\{(p_{i,j}(z_2,z_3))\  0\le i,j\le m_1 :\sum_{\matrix{i-j=k\cr|k|\le m_1} }
p_{i,j}(z_2, z_3)=p_k(z_2, z_3),\}.
$$
Note that the matrices in ${\cal F}$ are banded since $p_k =0,\ |k|>n$. We
look for a matrix ${\widehat P}(z_2,z_3)$ in ${\cal F}$ that is positive
definite for $|z_2|=1=|z_3|$.
The polynomial matrix $\widehat{P}(z_2,z_3)$ can be written as 
$$
\widehat{P}(z_2,z_3)=\sum_{k=-n}^{n}{\tilde P}_k(z_3) z_2^k,  
$$
where each ${\tilde P}_k(z_3)$ is an $(m_1+1)\times (m_1+1)$ Toeplitz matrix. Thus we
can write
$$
\widehat{P}(z_2,z_3)={\bf z_2}^\dagger{\bar P}(z_3){\bf z_2}, 
$$
where 
$$
{\bf z_2}=[I_{m_1},z_2 I_{m_1}, \ldots, z_2^{m_2} I_{m_1}]^T,
$$ 
with $I_{m_1}$ being the $(m_1+1)\times (m_1+1)$ identity matrix and $m_2\ge n$.
The polynomial ${\bar P}(z_3)$ is a matrix polynomial of size $(m_1+1)(m_2+1)\times (m_1+1)(m_2+1)$.
If it is nonnegative definite we can factor it into a polynomial matrix
$Q(z_2)$, i.e., ${\bar P}(z_3)=Q(z_3)^\dagger Q(z_3)$ by the matrix F\'ejer-Riesz theorem
(cf. [Helson'64] or [Mclean and Woerdeman'01]) then we have
$$
P(z_1,z_2,z_3)= \left(Q(z_3){\bf z_2}{\bf z_1}\right)^\dagger
\left(Q(z_3){\bf z_2}{\bf z_1}\right)
$$
which is a sum of square magnitudes of polynomials in $z_1, z_2, z_3$.

Our task then is to produce a positive definite polynomial matrix for any  
given positive multivariate Laurent polynomial. 
We resume our discussion on the two
variable case again and rewrite $P(z_1,z_2)$ as follows:
$$
P(z_1,z_2)= \sum_{k=-n_1}^{n_1}p_k(z_2)z_1^k={\bf z_{m_1}}^\dagger P_{m_1}(z_2) {\bf z_{m_1}} 
$$
where $m_1\ge n_1$, ${\bf z_{m_1}}=[1, z_1, z_1^2, \cdots, z_1^{m_1}]^T$, and 
$$
P_{m_1}(z_2)=[ p_{jk}(z_2)]_{0\le j, k\le m_1}
$$
with polynomial entries $p_{j,k}(z_2)$ given by
$$
p_{jk}(z_2)={1\over m_1+1-|j-k|} p_{k-j}(z_2), \forall j,k=0, \cdots, m_1. 
$$
Note that $p_{jk}(z_2)=0$ for $|j-k|>n_1$. 
Under this decomposition we can show that for some $m_1$ large enough, 
the matrix $P_1(z_2)$ will be positive definite when $P(z_1,z_2)$ is positive definite. 
To see this we note $P(z_1,z_2)$ is the symbol of 
the following bi-infinite Toeplitz matrix,
$$
\left[ \matrix{ \ddots & \ddots & \ddots & \ddots & \ddots &\ddots &\ddots  \cr
\ddots & p_0(z_2) & p_{-1}(z_2) & \cdots &p_{-n}(z_2) & 0 & \cdots \cr
\ddots & p_1(z_2) & p_0(z_2) & \ddots & \ddots & \ddots &\ddots \cr
\ddots & p_2(z_2) & p_1(z_2) & \ddots & \ddots & \ddots &\ddots \cr
\ddots & \ddots & \ddots & \ddots & \ddots & \ddots &\ddots\cr  
\ddots & p_{n}(z_2) & p_{n-1}(z_2) & \cdots & \ddots &\ddots &\ddots  \cr
\ddots & \ddots & \ddots & \ddots & \ddots &\ddots &\ddots  \cr}\right]. \eqno(4) $$
The positivity of $P(z_1,z_2)$ implies that any central section of the this matrix, i.e., any
square block with the diagonal consistent with the main diagonal 
$$
\hbox{diag}(\cdots, p_0(z_2), p_0(z_2), p_0(w_2), \cdots)
$$ 
is positive as explained at the beginning of this section.  Typically, we have
$$
p_0(z_2)>0, \quad \left[\matrix{p_0(z_2) & p_{-1}(z_2)\cr p_1(z_2) & p_0(z_2)\cr}\right]>0, 
\quad \left[\matrix{ p_0(z_2) & p_{-1}(z_2) & p_{-2}(z_2)\cr p_1(z_2) & p_0(z_2) & p_{-1}(z_2)\cr
p_2(z_2) & p_1(z_2) & p_0(z_2)\cr}\right]>0, \cdots .
$$
For convenience, we denote by ${\cal P}_2$ and ${\cal P}_3$ to be the $2\times 2$ and 
$3\times 3$ matrices above, respectively. In general, we use ${\cal P}_k$ to denote the
$k\times k$ central block matrix from the bi-infinite Toeplitz matrix (4) above. 

We now look at the matrix $P_{m_1}(z_2)$ given by,
$$
\left[ \matrix{{1\over {m_1}+1} p_0(z_2) & {1\over {m_1}}p_{-1}(z_2) & \cdots & {1\over {m_1}+1-n_1}
p_{-n_1}(z_2) & 0 & \cdots \cr
{1\over {m_1}}p_1(z_2) & {1\over {m_1}+1}p_0(z_2) & {1\over {m_1}}p_{-1}(z_2) & \ddots & \ddots &\ddots \cr
{1\over {m_1}-1}p_2(z_2) & {1\over {m_1}}p_1(z_2) & \ddots & \ddots & \ddots &\ddots \cr
\vdots & \ddots & \ddots & \ddots & \ddots &\ddots\cr  
{1\over {m_1}+1-n_1}p_{n_1}(z_2) & \ddots & \cdots & \ddots &\ddots &\ddots  \cr
0 & \ddots & \ddots & \ddots &\ddots &\ddots  \cr
\vdots & \ddots & \ddots & \ddots & \ddots &{1\over {m_1}+1}p_0(z_2)  \cr}\right].
$$
Note that each diagonal sums to $p_{i-j}(z_2)$ so $P_{m_1}$ is in
${\cal F}(z_2)$ where ${\cal F}(z_2)$ is defined as above with the obvious modifications.
With $x=[x_0,x_1,\cdots, x_{m_1}]^T, $ 
we need to prove that $x^* P_{m_1}(z_2)x>0$. First we write 
$$
x^\dagger P_{m_1}(z_2)x = {1\over {m_1}+1} x^\dagger {\cal P}_{m_1} x + 
{1\over {m_1}+1}x^\dagger R_{m_1} x
$$ 
with a remainder matrix $R_{m_1}$. The $\ell_2$ norm of $R_{m_1}$ can be estimated
using the column norm to give
$$
\|R_{m_1} \|_2 \le  {n_1(n_1+1)C_1\over2(m_1-n_1)},
$$
where $C_1=\sup_{i,|z_2|=1}|p_i(z_2)|$. If $P(z_1,z_2)\ge \epsilon$ then
$x^\dagger{\cal P}_{m_1} x\ge \epsilon\|x\|_2$, 
so that if $\displaystyle { n_1(n_1+1) C\over2(m_1-n_1)} 
<\epsilon$, then $x^\dagger
P_{m_1}(z_2)x>0$. Then an application of the matrix Riesz-Fejer
Theorem  yields

\proclaim Theorem 2.1. Let $P(z_1,z_2)=\displaystyle \sum_{k=-n_1}^{n_1} p_k(z_2)z_1^k\ge \epsilon>0$ 
be strictly positive on bi-torus $|z|=1=|z_2|$. Then $P(z_1,z_2)$ can be factored into a sum of
squares of polynomials in $z_1$ and $z_2$. 
The total number of terms in the sum is less than or equal to $m_1+1$ 
with $m_1$ being  an integer such that 
$$
{n_1(n_1+1) C_1\over2(m_1-n_1)}< \epsilon,
$$
and the degrees of each of the polynomials is bounded by $m_1$ in $z_1$ and $n_2$ in $z_2$.

We remark that when $P(z_1,z_2)$ has different coordinate degrees $n_1, n_2$, 
it may be worthwhile depending upon $C_1$ to choose the smaller among $n_1$ and $n_2$ 
in order to have a fewer terms in the sum of square magnitudes of polynomials for $P(z_1,z_2)$.

Next we generalize the result in Theorem 2.1 to  
the multivariate setting which is known from [Dritschel'04]. 

\proclaim Theorem 2.2(Dritschel). Let $P(z_1, \cdots, z_d)$ be a
multivariate  Laurent polynomial which 
is strictly positive on the multivariate torus $|z_1|=|z_2|=\cdots=|z_d|=1$, where $d\ge 2$ is 
an integer. Then $P(z,w)$ can be expressed as a sum of square magnitudes of 
polynomials in $z_1, \cdots, z_d$.
 
\pf We shall use the arguments in the proof 
of the previous Theorem. Write $P(z_1,z_2,\ldots,z_d)=P(z_1,z)=\sum_{j=-n_1}^{n_1}
p_j(z)z_{1}^j>0$, where $z$ is the usual multivariable notation beginning
with $z_2$. We know that $P(z_1,z)$ is
the symbol of the bi-infinite Toeplitz matrix given by (4) with $z_2$
replaced by the multivariable $z$.
It follows that any central section along the main diagonal is strictly positive definite 
as explained before. Write 
$$
P(z_1,z)= {\bf z_1}^\dagger P_{m_1}(z){\bf z_1}, \eqno(5)
$$
where ${\bf z_1}$ given by equation (3) and 
$P_{m_1}(z)=\left[p_{j,k}\right]_{0\le j, k\le m}$ 
is a matrix of size $(m_1+1)\times (m_1+1)$ with entries 
$$
p_{jk}={1\over m_1+1-|j-k|} p_{j-k}(z), \quad \forall j, k=0, 1, \cdots, m_1.
$$
If $P>\epsilon$ the argument in Theorem 2.1 shows that for $m_1$ large
enough there is an $\epsilon_1>0$ such that $x^\dagger P_{m_1}(z)x>\epsilon_1$ on
the d-1 torus if $\displaystyle {n_1(n_1+1){\hat C}_1\over2{m_1-n_1}}<\epsilon$, where in
this case ${\hat C}_1 = \sup_{i, |z_j|=1, j=2,..., d}|p_i(z)|$. Write
$\displaystyle P_{m_1}(z_2,z')=\sum_{k=-n_2}^{n_2}{\tilde p}_k(z')z_2^k$, where ${\tilde p}_k$ are
$(m_1+1)\times (m_1+1)$ Toeplitz matrices and $z'=(z_3,\ldots, z_d)$. Now set
$$
\hat{p}_{jk}={1\over m_2+1-|j-k|} \tilde{p}_{j-k}(z'), \quad \forall j, k=0, \cdots, m_2
$$
with $m_2\ge m_1$ and $\displaystyle P_{m_2}(z')=\left[ \hat{p}_{j,k} \right]_{0\le j, k\le m_2}$.
As above we have that
$$
x^\dagger P_{m_2}x={1\over m_2+1}x^\dagger{\cal P}_{m_2} x+{1\over m_2+1}x^\dagger{\cal R}_{m_2}x.
$$
As above the norm of ${\cal R}_{m_2}$ can be bounded by
$$
\|{\cal R}_{m_2}\|_{2}\le {n_2(n_2+1)C_2 \over2(m_2-n_2)},
$$
where $\displaystyle C_2=\sup_{i,|z_2|=\cdots |z_d|=1}\|{\tilde p}_i(z')\|_2$. Thus for
$m_2$ sufficiently large, $P_{m_2}$ is a positive matrix polynomial. We 
continue the process until we arrive at the positive trignometric matix
polynomial $P_{m_{d-1}}(z_d)$ which can be factored by the matrix F\'ejer-Reisz
Theorem. We have thus established  the proof. \eop

Note that the number of factors will be $(m_1+1)(m_2+1)\cdots (m_{d-1}+1)$
and the degrees of the polynomials at most $m_1$ for $z_1$ ... $m_{d-1}$
for $z_{d-1}$ and $n_d$ for $z_d$ . We  note that we could have avoided the
use of the matrix Fej\'er-Riesz lemma by eliminating all variables then using a
square root of a positive matrix (see Mclean and Woerdeman'01). We will
consider an alternative computationally attractive method for computing
factorizations in the next section.

\sect{3}{Computing Approximate Factorizations}
As shown in the previous section, an important step in the factorization of 
multivariate Laurent polynomials is to compute the factorization of univariate polynomial matrices. 
 Recall  a computational
algorithm for factorizations of one variable trignometric
polynomials was exploited in [Lai'94]. This method can be developed
to factorize polynomial matrices in the univariate setting.   
Let us first introduce  some necessary notation and definitions in order to explain the method
in more detail. 

Let $\ell_2$ stand for the space 
of all square summable sequences. Let $\|{\bf x}\|_2$ denote the standard norm on $\ell_2$.  
We note that  for any operator $A$ from $\ell^2 \mapsto \ell^2$, 
$A$ can be expressed by a bi-infinite matrix. 

\proclaim Definition 3.1.  A bi-infinite matrix $A=(a_{ik})_{i,k\in {\bf Z}}$ is said to be of 
exponential decay off its diagonal  if 
$$
\|a_{ik}\|_2 \le K r^{|i-k|}
$$
for some constant $K$ and $r\in (0,1)$, where ${\bf Z}$ is the collection of 
all integers.   $A$ is banded with band width $b$ if $a_{ik}=0$ for all 
$i, k\in {\bf Z}$ with $|i-k|> b$.

If $A$ is a positive operator, then  
there exists the unique  positive bi-infinite square root matrix $Q$ of $A$ such that $Q^2=A$. 
If $A=B^\dagger B$ for another bi-infinite matrix $B$, then 
there exists a unitary matrix $U$ such that $B=UQ$. 

Recall from the previous section that  given any Laurent polynomial $P(z)$, 
we can view $P(z)$ to be the symbol of  a bi-infinite Toeplitz matrix ${\cal P}$. 
The computational scheme introduced in [Lai'94] roughly speaking is to choose a central section 
$$P_N=(p_{j-k})_{-N\le j,k\le N}$$  
of matrix ${\cal P}$ and compute a Cholesky factorization i.e
$P_N=C_N^\dagger C_N $ where $C_N$ is an upper triangular matrix with
positive diagonal entries,
if $P_N$ is positive definite or use the singular value decomposition (SVD) to 
find $Q_N$ if $P_N$ is nonnegative definite and then find a Householder matrix $H_N$ such that
$C_N=H_NQ_N$ is upper triangular. Then the nonzero entries in the middle row of $C_N$   
approximate that in the middle row (in fact any row) of ${\cal C}$ whose symbol $C(z)$ is a 
factorization of $P(z)$, i.e., $P(z)=C(z)^* C(z)$. 

For the extension of this method to matrix polynomials, let $$\ell^m_k=\{{\bf
x}=\{x_i\}_{i\in {\bf Z}}, x_i\in {\bf R}^m, \|{\bf x}\|_k <\infty\},\quad  k=1,2$$
and $B(\ell^m_2)$ be the set of bounded linear operators on $\ell^m_2$. Let
$\Pi_N\in B(\ell^m_2)$ be the projection given by
$$
\Pi_N{\bf x}={\bf y},\ {\bf y}=\{y_i\}: y_i=0, |i|>N, y_i=x_i,\ |i|\le N.
$$
If $P\in B(\ell^m_2)$ is positive definite we will be interested in considering the 
$(2N+1)m\times (2N+1)m$ submatrix of $P$ centered at the index zero which will be called 
the Nth central section and 
which is also positive definite.  We will also be interested in extensions of
various finite matrices $A_N$ to $B(\ell^m_2)$ given by
$$
\left[\matrix{0&0&0\cr0& A_N&0\cr 0&0&0}\right],
$$
which with a slight abuse of notation will also be called $A_N$.

Consider the matrix polynomial $P(z)= \sum_{j=-n}^n p_{j}z^j$
with matrix coefficients $p_k$'s of size $m\times m$, 
then ${\cal P}=(p_{i-j})_{i,j\in {\bf Z}}\in B(\ell^m_2)$ defined by $m\times m$ 
matrix blocks $p_k, -n\le k\le n$ is a bi-infinite block Toeplitz matrix 
whose symbol is $P(z)$.
 As shown earlier if $P(z)$ is Hermitian  nonnegative definite, so is 
${\cal P}$. 
Let $C(z)$ be a factorization  of $P(z)$  i.e., $P(z)=C (z)^\dagger$ C(z), 
then ${\cal P}={\cal C}^\dagger{\cal C}$, where ${\cal C}$ is a bi-infinite upper triangular banded 
block Toeplitz matrix associated with $C(z)$.   
On the other hand, if ${\cal P}={\cal C}^\dagger{\cal C}$ for a upper
triangular banded block Toeplitz matrix, 
then the symbol $C(z)$ of ${\cal C}$ satisfies $P(z)=C(z)^\dagger C(z)$. If
$P(z)$ is positive definite then it follows from the matrix Fej\'er-Riesz 
Lemma [Helson'64], [Mclean-Woerdeman'01] that it is possible to 
choose ${\cal C}$ so that it has positive diagonal entries. We shall prove
the following,
\proclaim Theorem 3.1. Let $P(z)=\sum_{-n}^n p_k z^k$ be an $m\times m$ 
matrix polynomial that is positive definite for $|z|=1$. Let ${\cal
P}=(p_{i-j})_{i,j,\in {\bf Z}}={\cal C}^\dagger{\cal C}$ where ${\cal C}$
is an upper triangular banded block Toeplitz with positive diagonal entries, 
$P_N$ be the $Nth$ central section 
of ${\cal P}$, and $C_N$ its  Cholesky factor (which we extend as described
above). Then
$$
\|(\hat C_N-{\cal C}_N)\delta\|_2<K\rho^N,
$$
for some $\rho\in (0,1)$, where  $\delta\in \ell^m_2$ is a vector with a finite number
of nonzero entries.

\noindent
{\bf Remark} For the numerical computation below we will choose $\delta$
with zero components  except for $\delta_0=I_m$, the $m\times m$
identity matrix.

The proof of  Theorem 3.1 is based upon the following Theorem 3.2 and Lemmas 3.3 and 3.4.

\proclaim Theorem 3.2. Suppose that $A\in B(\ell_2)$ is a positive banded operator such
that $\Vert A-I\Vert_2<1$. Let $Q$ be the unique positive square root of $A$, $A_N$ be a central section,
and ${\hat Q}_N$ be the positive  matrix such that  ${\hat Q}_N^2
=A_N$. Then
$$
\Vert  (Q-{\hat Q}_N)\delta \Vert_2 \le K \lambda^{N} \eqno(6)
$$
for some $\lambda\in (0,1)$ and a positive constant $K$. In equation (6) $\delta \in
l_2$ is any vector with a fixed number of nonzero entries.   

To prove the above Theorem, we begin with the following lemmas. 

\proclaim Lemma 3.3.  Suppose that $A$ is banded with bandwith $b$ and $\Vert A-I\Vert_2\le r<1$. 
Then $Q=(q_{ik})_{i,k\in {\bf Z}}$ then $|q_{l,k}|\le K r^{{|l-k|\over
b}}$. If $A$ is invertible, then the entires of $Q^{-1}$ satisfy a similar bound.  

\pf We only prove the exponential decay property of $Q$. The proof of that
of $Q^{-1}$ is similar. 
The uniqueness of $Q$ and the convergence of the following series 
$$
 \sum_{i=0}^\infty (-1)^i {(2i-3)!!\over (2i)!!} (A-I)^i
$$
implies that 
$$
Q= \sqrt{A}=\sqrt{I+(A-I)} =
 \sum_{i=0}^\infty (-1)^i {(2i-3)!!\over (2i)!!} (A-I)^i.
$$
$A$ is banded and so is $A-I$. If $A-I$ has bandwidth $b$, then 
$(A-I)^i$ is also banded with  bandwidth $ib$. Thus, 
$$
q_{jk}= \sum_{i\ge |j-k|/b}^\infty (-1)^i {(2i-3)!!\over (2i)!!} (A-I)^i_{jk}, 
$$
where $(A-I)_{jk}$ denotes the $(j,k)^{\rm th}$ entry of $A-I$ and similar for $(A-I)^i_{jk}$. 
It follows that 
$$
|q_{jk}|  \le  K r^{|j-k|/b}
$$
for some constant $K$. This finishes the proof. \eop

Let us write 
$$
Q= \left[\matrix{ \alpha_1& B & \alpha_2 \cr B^\dagger  & Q_N & C^\dagger  \cr \alpha_3 & C &
\alpha_4\cr}\right] \hbox{  and  }
A=\left[\matrix{ \beta_1& a & \beta_2 \cr a^\dagger  & A_N & c^\dagger  \cr \beta_3 & c &
\beta_4\cr}\right].
$$
Note that  $Q^2=A$ implies $A_N= Q_N^2+B^\dagger  B+C^\dagger  C$  or ${\hat Q}_N^2-Q^2_N=
B^\dagger  B+C^\dagger  C$ where ${\hat Q}_N^2=A_N$. Thus, we have 
$$
(Q_N+{\hat Q}_N)({\hat Q}_N-Q_N)= {\hat Q}_N^{2}-Q_N^2+Q_N{\hat Q}_N-{\hat Q}_N
Q_N=B^\dagger  B+C^\dagger  C+R, \eqno(7)
$$
where $R$ is defined in the following, 

\proclaim Lemma 3.4. (cf. [Lai'94])  Let $R=(r_{jk})_{-N\le j,k\le N}:=Q_N{\hat Q}_N-{\hat Q}_N Q_N$. 
Then $r_{jk}=O(r^{N/(4b)})$ for $k=-N/4+1,\cdots, N/4-1$ and $j=-N,\cdots,N$. 
 
\pf (of Theorem 3.2.) From equation (6) we find that, $({\hat Q}_N-Q_N)=(Q_N+{\hat Q}_N)^{-1}(B^\dagger  B+C^\dagger  C+R)$. By Lemma 3.3., we 
can prove that the entries of $B^\dagger  B+C^\dagger  C$ have the exponential decay property: 
$(B^\dagger  B+C^\dagger  C)_{jk} =O(r^{N-|k|}),-N\le k\le N$.

The positivity of $A$ implies that $Q$ is positive and so is $Q_N$. It follows that $\|Q_N^{-1}\|_2$
is uniformly bounded, furthermore since $hat Q_N$ ia also nonegative we find,
$$
\| (Q_N + \hat{Q}_N)^{-1}\|_2 \le \|Q_N^{-1}\|_2 \le K_1 <\infty
$$
for a positive constant $K_1$ independent of $N$, where we have used the fact that
$\hat{Q}_N$ is nonnegative.  Therefore, we conclude that 
$$\eqalign{
\|({\hat Q}_N-Q_N)\delta_N\|_2 \le & 
\|(Q_N+{\hat Q}_N)^{-1}\| \|(B^\dagger B+C^\dagger C+R)\delta_N\|_2 \cr 
\le & K_1 \|(B^\dagger B+C^\dagger C+R)\delta_N\|_2 \cr}
$$
where $\delta_N$ is the finite vector whose entries match those of $\delta$. The proof is  completed by
extending $Q_N,\ {\hat Q_N}$, replacing $\delta_N$ by $\delta$, and
noticing that by Lemma 3.3 $\|(Q_N-Q)\delta\|_2<K_1 \lambda^N$, $\lambda<1$. \eop  

\pf (of Theorem 3.1) Suppose that 
$$\sup_{|z|=1}\|P(z)\|_2<1. \eqno(8)$$ Otherwise divide $P$
by a sufficiently large constant so that (8) holds. Let $Q$ be the unique positive square root of
${\cal P}$, and $Q_N$ the positive square root of
$P_N$. From Theorem 3.2 we know that $\|(Q_N-Q)\delta\|_2<K\rho^N$ with
$\rho<1$. Let $U$ be the unitary matrix such that ${\cal C}=UQ$. Then
$$
\|(Q_N-Q)\delta\|_2=\|(UQ_N-{\cal C})\delta\|_2.
$$
Write $UQ_N=\tilde Q_N+ L^1_N$ where $\tilde Q_N$ is upper
triangular and $L^1_N$ is strictly lower triangular, then $UQ_N =q_N
+l_N$ where $q_N=\Pi_N \tilde Q_N\Pi_N^\dagger$ and $l_N=L^1_N+{\tilde Q}_N-q_N$. 
Theorem 3.2 shows that $\|l_N\delta\|_2$ tends
to zero exponentially fast. Furthermore  since ${\tilde Q}_N$ is symmetric, 
$$\eqalign{
P_N=&{\tilde Q}_N^2={\tilde Q}_N^\dagger {\tilde Q}_N = (U{\tilde Q}_N)^\dagger (U {\tilde Q}_N)\cr
=&(q_N+l_N)^\dagger (q_N+l_N)\cr
=& q_N^\dagger q_N + l_N^\dagger q_N +q_N^\dagger l_N +l_N^\dagger l_N. \cr}
$$
That is, we have 
$$
C_N^\dagger C_N - q_N^\dagger q_N = l_N^\dagger q_N +q_N^\dagger l_N +l_N^\dagger u_N.
$$
Since $Q_N$ is uniformly bounded so is $q_N$ and we find,
$$
\|(C_N^\dagger C_N -q_N^\dagger q_N)\delta\|_2<K_2\lambda^N.
$$
Restricting the above quantities to their finite matrices we
note because of the strict positivity of $P$, $\|C_N\|_2$ is uniformly
bounded from below hence $C_N^{-1}$ is uniformly bounded. Furthermore since
$C_N$ has the same size as $q_N$,
$$
\|(I -(C_N^\dagger)^{-1}q_N^\dagger q_N C_N^{-1})\delta_N\|_2<K_3\lambda^N,
$$
where $\delta_N=C_N\delta$ for any $\delta$ with finitely many nonzero entries. 
Note that the factor $q_N C_N^{-1}$ is upper triangular while $(C_N^\dagger)^{-1}q_N^\dagger$ is   
lower triangular.  The above inequality shows that
$\|(q_N C_N^{-1}-I)\delta\|_2<K_3\lambda^N$. This completes the proof. \eop

\sect{4}{Numerical Examples}

In this section we give three examples to illustrate how the computational method works for 
polynomial matrix factorizations.

\noindent
{\bf Example 4.1.}  We first consider a univariate polynomial matrix
$$
P(z):=\left[ \matrix{ 8+z +1/z & 1+z\cr 1+1/z & 1\cr}\right].
$$
It is clear that the matrix is Hermitian and positive definite. We write 
$$
P(z)=\left[\matrix{ 8& 1\cr 1 & 1\cr}\right]+ \left[\matrix{1&1\cr 0& 0\cr}\right]z 
+\left[\matrix{1&0\cr 1& 0\cr}\right]/z.
$$
We assemble a bi-infinite Toeplitz matrix whose $10\times 10$ block is as shown below.
$$
\left[\matrix{8& 1& 1& 1& 0& 0& 0& 0& 0& 0\cr
1& 1& 0& 0& 0& 0& 0& 0& 0& 0\cr
1& 0& 8& 1& 1& 1& 0& 0& 0& 0\cr
1& 0& 1& 1& 0& 0& 0& 0& 0& 0\cr
0& 0& 1& 0& 8& 1& 1& 1& 0& 0\cr
0& 0& 1& 0& 1& 1& 0& 0& 0& 0\cr
0& 0& 0& 0& 1& 0& 8& 1& 1& 1\cr
0& 0& 0& 0& 1& 0& 1& 1& 0& 0\cr
0& 0& 0& 0& 0& 0& 1& 0& 8& 1\cr
0& 0& 0& 0& 0& 0& 1& 0& 1& 1\cr}\right].
$$
We use the Cholesky decomposition method to a $20\times 20$ central block and get a lower triangular
matrix $F$. Let $P0$ be  the right and bottom $2\times 2$ block from $F$ which is
$$
P0 := \left[\matrix{{\sqrt{385}\over 7}&0 \cr {6\over \sqrt{385}}
&{\sqrt{2310}\over 55}\cr}\right].
$$
Choose the $2\times 2$ block  next to $P0$ as follows
$$
P1 := \left[\matrix{{\sqrt{385}\over 55}&{-\sqrt{2310}\over 385}\cr{\sqrt{385}\over 55}& {-\sqrt{2310}\over 385}\cr}\right]  
$$
Define $Q^\dagger(z)=P0+P1/z$ and then we have $P(z)=Q(z)^\dagger Q(z)$.  \eop 

\noindent
{\bf Example 4.2.} We next consider a bivariate polynomial 
$$
\eqalign{
&P(x,y)= 41+5x^2+5y^2+15/x+20/y+5/x^2+5/y^2+15x+20y+5xy\cr
&\qquad +8y/x+5/(xy)+8x/y+2x/y^2+3y/x^2+3x^2/y+x^2/y^2+2y^2/x+y^2/x^2\cr}
$$
It is a positive polynomial since $P(x,y)=p(x,y)p(1/x,1/y)$ with 
$p(x,y)=5+2x+3y+xy+x^2+y^2$. Let us write 
$$
P(x,y)=[1,1/x,1/x^2] \widetilde{P}(y)\left[\matrix{1\cr x\cr x^2\cr}\right],
$$
with 
$$
\eqalign{&\widetilde{P}(y):= \cr
&\left[\matrix{
{41\over 3}+{5y^2\over 3}+{20\over 3y}+
{5\over 3y^2}+{20\over 3}y & {15\over 2}+{5\over 2}y+{4\over y}+{1\over y^2}& 
5+{3\over y}+{1\over y^2}\cr
{15\over 2}+4y+{5\over 2y}+y^2& {41\over 3}+{5\over 3}y^2+{20\over 3y}+{5\over 3y^2}+{20\over 3}y
& {15\over 2}+{5\over 2}y+{4\over y}+{1\over y^2}\cr 
5+3y+y^2 & {15\over 2}+4y+{5\over 2y}+y^2& {41\over 3}+{5\over 3}y^2+{20\over 3y}+
{5\over 3y^2}+{20\over 3}y\cr}\right]. \cr}
$$
Then we can write 
$$
P(x,y)= [1,1/x,1/x^2] \widetilde{P}(y) \left[\matrix{1\cr x\cr x^2\cr}\right],
$$
where $\widetilde{P}(y)=\sum_{j=-2}^2 p_j y^j$ with $p_{-2}, \cdots, p_2$ being given below:
$$\eqalign{
p_0=&\left[\matrix{{41\over 3}& {15\over 2}& 5\cr {15\over 2} & {41\over 3} & {15\over 2}\cr
5 & {15\over 2} & { 41\over 3}}\right], 
p_1=\left[\matrix{{20\over 3} &{ 5\over 2} & 0\cr
4 & {20\over 3} & {5\over 2}\cr 
3 & 4 & {20\over 3}}\right],  p_{-1}=p_1^\dagger ,\cr 
p_2=&\left[\matrix{{5\over 3} & 0& 0\cr 1 & {5\over 3} & 0\cr 1& 1& {5\over 3}\cr}\right], 
p_{-2}=p_2^\dagger .\cr}
$$
We now assemble a bi-infinite Toeplitz matrix whose $9\times 9$ central block are shown as follows:
$$
\left[\matrix{{41\over 3}& {15\over 2} & 5 & {20\over 3} &{5\over 2}& 0&  {5\over 3} & 0& 0\cro
{15\over 2} & {41\over 3} & {15\over 2} & 4 & {20\over 3} & {5\over 2} & 1 & {5\over 3} & 0\cro
5 & {15\over 2} & {41\over 3} & 3 & 4 & {20\over 3} & 1 & 1& {5\over 3}\cro
{20\over 3} & 4 & 3 & {41\over 3} & {15\over 2} & 5 & {20\over 3} & {5\over 2} & 0\cro
{5\over 2} & {20\over 3} & 4 & {15\over 2} & { 41\over 3} &{15\over 2} &4 &{20\over 3}& {5\over 2}\cro
0 & {5\over 2} & {20\over 3} & 5& {15\over 2} & {41\over 3} & 3& 4& {20\over 3}\cro
{5\over 3} & 1 & 1& {20\over 3} & 4& 3& {41\over 3} &  {15\over 2} & 5\cro
0 & {5\over 3} & 1& {5\over 2} & {20\over 3} & 4& {15\over 2} & {41\over 3} & {15\over 2}\cro
0 & 0& {5\over 3} & 0&  {5\over 2} & {20\over 3} & 5& { 15\over 2} & {41\over 3}}\right].
$$
We use the Cholesky factorization of a central block matrix of size $60\times 60$. Let $F$ be the
lower triangular factorization. Then choose $Q_0$ to 
be the $3\times 3$ block at the bottom and right of $F$, $Q_1$ the $3\times 3$ block next to $Q_1$
and $Q_2$ the $3\times 3$ block next to $Q_1$ that is 
$$
\eqalign{
Q_0=& \left[\matrix{3.185602126 & 0 & 0\cr
1.873651218& 2.539725049& 0\cr
1.524622962& 1.128505745& 2.269126602\cr}\right], \cr
Q_1=& \left[\matrix{1.797364251& 0.08381502303 & -0.0003518239229\cr
0.7675275947& 1.633796832& 0.06150315980\cr
0.00008111923034& 0.9665117592& 1.856367398\cr}\right]\cr
Q_2=& \left[\matrix{0.5231873284& 0.007768330871 & 0.08530594055\cr
0 & 0.6562390159 & 0.1143305535\cr
0 & 0 & 0.7344969935\cr}\right]. }
$$
Let $Q(y)^\dagger=Q_0+Q_1/y+Q_2/y^2$ and then $Q(y)^\dagger Q(y) \approx \widetilde{P}(y)$. In fact 
the maximum error of each entry of $Q(y)^\dagger Q(y)- \widetilde{P}(y)$ is less than or equal to $10^{-8}$.
\eop 

\noindent
{\bf Example 4.3.} Let us consider a bivariate polynomial which has a zero on the bi-torus:
$$
P(x,y)= 30+14/x+11/y+4/x/y+14x+6x/y+11y+6y/x+4xy.
$$
It is the product of $P(x,y)=(4+3x+2y+1)(4+3/x+2/y+1)$ which is zero at $x=-1, y=-1$. We write
$$
P(x,y)= p_0(y)+p_1(y)x+p_{-1}(y)/x
$$
for $p_0(y)=30+11/y+11y$, $p_1(y)=14+6y+4/y$, and $p_{-1}(y)=14+4y+6/y$. It is the symbol of 
an bi-infinite Toeplitz matrix. One of its central section is as shown below.
$$
\left[\matrix{11/y+30+11y & 4/y+14+6y & 0 & 0\cr
6/y+14+4y & 11/y+30+11y & 4/y+14+6y& 0\cr
0& 6/y+14+4y& 11/y+30+11y & 4/y+14+6y\cr
0& 0& 6/y+14+4y& 11/y+30+11y\cr}\right].
$$
Since $P(x,y)$ has no simple factors (see the next section),  any central sections of  
the bi-infinite Toeplitz matrix is positive by Lemma 5.1. We consider 
several central sections $P_m$ of size $m=16\times 16$, $32\times 32$, $64\times 64$ and 
$128\times 128$.
For each of these central sections, $P_m$ is a univariate polynomial in $y$ with matrix coefficients
and $P_m(y)$ is positive. Thus, $P_m(y)=Q_m(y)^\dagger Q_m(y)$. To compute $Q_m(y)$, we use
the computational method in \S 3 to yield an approximation $\tilde{Q}_m$ of $Q_m$. 
As the size of central sections increases, the
$Q_m$ converges to the corresponding entries in the bi-infinite Toeplitz 
matrix.  We use the entries 
on the last row of $\tilde{Q}_m$ to construct an approximation of $Q_m(y)$ and hence the factorization
of $P(x,y)$ and listed below. 
$$
\left[\matrix{ {\rm size} & {\rm factorization} \cr
16\times 16 &4.01207952+2.984741799x+2.000226870y+0.996712925xy\cr
32\times 32 &4.004041536+2.994924757x+2.000034879y+0.998949058xy\cr
64\times 64 &4.001381387+2.998269650x+2.000005690y+0.999648058xy\cr
128\times 128 &4.00069369+2.999134582x +1.99999896y+0.999821915xy\cr}\right].
$$
As we know that the factorization is $4+3x+2y+1$,  the approximations are very good. \eop

\sect{5}{Nonegative bivariate Trignometric Polynomials}
Finally we consider the problem of factorization of nonnegative 
multivariate polynomials.  Let us start with $P(z, w)\ge 0$. If for some $z_0$ with $|z_0|=1$, 
$P(z_0,w)=0$ for all $w$ with $|w|=1$, we say that $P(z,w)$ has a simple factor at $z_0$. If 
$P(z,w)$ has a simple factor at $z_0$, then $P(z,w)$ has factors $(z-z_0)$ and $(1/z-1/z_0)$. Let
us factor them out. Then $P(z,w)/((z-z_0)(1/z-1/z_0))$ is still nonnegative. Similarly, if 
$P(z,w_0)=0$ for all $z$ with $|z|=1$, $P(z,w)$ has a simple factor at $w_0$. In this case,
$P(z,w)$ has two factors $(w-w_0)$ and  $(1/w-1/w_0)$ which can be factored out from $P(z,w)$. 
 Without loss of generality, we may assume that $P(z,w)\ge 0$ does not have
any simple factors. Writing $\displaystyle P(z,w)= \sum_{j=-n}^n p_j(w)z^j$, we view that $P(z,w)$ is a 
polynomial of $z$ and it is the symobl of a bi-infinite Toeplitz matrix in (4) with $w$ in place
of $z_2$. We have the following 

\proclaim Lemma 5.1. Suppose that $P(z,w)\ge 0$ does not have any simple factors. Then any central 
section of the bi-infinite Toeplitz matrix in (4) is strictly positive definite. 

\pf Since $P(z,w)\ge 0$, we know that any central section of the matrix in (4) is nonnegative
definite. Suppose that a central section $T_m(w)$ of the matrix in (4) is not positive definite
for  $w=w_0$. 
Then  there exists a vector ${\bf x}$ such that $T_m(w_0){\bf x}=0$, i.e.,  
${\bf x}^\dagger T_m(w_0) {\bf x}=0$. Thus, we have, for $z=e^{i\theta}$, 
$$
0= {\bf x}^\dagger T_m(w_0) {\bf x} = {1\over 2\pi}\int_0^{2\pi}F({\bf x})^*P(z,w_0)F({\bf x})d\theta. 
$$
It follows that
$$
|F({\bf x})|^2 P(z,w_0) =0, \quad a.e. 
$$
and hence, $P(z,w_0)\equiv 0$ since $|F({\bf x})|\not=0, a.e.$ and $P(z,w_0)$ is a Laurent polynomial.
That is, $P(z,w)$ has a simple factor at $w_0$. This contradicts the assumption on $P(z,w)$. \eop

 
Thus, for a central section $P_m$  of size 
$m\times m$ in the matrix in (4),  $P_m$ is positive. Since $P_m$ is a matrix polynomial in $w$,
by the matrix F\'ejer-Riesz factorization theorem (cf. [Helson'64]), 
$P_m$ can be factorized into $Q_m$, i.e., $P_m(w)=Q_m(w)^\dagger Q_m(w)$. Intuitively,  
the polynomial $Q_m$ is a good approximation of the factorization of the bi-infinite Toeplitz 
matrix ${\cal P}$ in (4) as $m$ sufficiently large. 
In the previous section, we presented an example (Example 4.3.) of $P(z,w)$ which is 
nonnegative without simple factors. 
Using our symbol approximation method, we compute an approximation of the factorization of 
$P_m$ for $m=16, 32, 64,$ and $128$. 
The numerical computation shows the factorizations converge. 


Let us now discuss the convergence a little bit more in detail. 
For simplicity, let ${\cal A}$ be a bi-infinite Toeplitz matrix associated with a univariate Laurent
polynomial $A(z)$ and ${\cal A}_N=(a_{jk})_{-N\le j, k\le N}$ be a central section of size 
$(2N+1)\times (2N+1)$ for a positive integer $N$.  
Suppose that each ${\cal A}_N$ is strictly positive. Thus we can 
obtain a factorization ${\cal A}_N={\cal B}_N^*{\cal B}_N$ by Cholesky factorization.  

\proclaim Lemma 5.2. For any
${\bf x, y}\in \ell_2$, ${\bf x_N}^\dagger {\cal A}_N {\bf y}:={\bf x}_N^\dagger{\cal A}_N {\bf y}_N$ 
converges to ${\bf x}^\dagger {\cal A}{\bf y}$ as $N\longrightarrow +\infty$, 
where ${\bf x}_N=(x_{-N}, \cdots, x_0, \cdots, x_N)^\dagger $
is the central section of size $2N+1$ of ${\bf x}$ around the index $0$ and similar for $y_N$.

\pf  For an integer $N>0$, 
$$\eqalign{
& {\bf x}^\dagger {\cal A}_{N} {\bf y} - {\bf x}^\dagger {\cal A} {\bf y} \cr
=& {1\over 2\pi}\int_0^{2\pi} \left(F({\bf x}_N)^* A(z) F({\bf y}_N)-
F({\bf x})^* A(z) F({\bf y})\right)d\theta \cr
=& {1\over 2\pi}\int_0^{2\pi} \left(F({\bf x}_N)-F({\bf x})\right)^* A(z) F({\bf y}_N)d\theta \cr
& + {1\over 2\pi}\int_0^{2\pi} F({\bf x})^* A(z)\left(F({\bf y}_N)- F({\bf y})\right)d\theta\cr}
$$
where $z=e^{i\theta}$. In the first inequality we used the fact that
${\bf x}^\dagger{\cal A}_N {\bf x}=(\Pi_N{\bf x})^\dagger{\cal A}\Pi_N{\bf
x}$ where $\Pi_N$ is the projection defined in section 3. Thus
$$
\eqalign{
&|{\bf x}^\dagger {\cal A}_{N} {\bf y} - {\bf x}^\dagger {\cal A} {\bf y} |\cr 
\le&  \|{\bf x}-{\bf x}_N\|_2 \|A(z)\|_\infty \|{\bf y}\|_2+ 
\|{\bf y}-{\bf y}_N\|_2 \|A(z)\|_\infty \|{\bf x}\|_2\cr
\longrightarrow &0\cr}
$$
as $N\rightarrow +\infty$. Here, $\|A(z)\|_\infty$ denotes the maximum norm of 
$A(z)$ over the circle $|z|=1$. This completes the proof. \eop

A consequence of the above Lemma 5.2 is that $\|{\cal B}_N {\bf x}\|^2_2$ converges to ${\bf x}^\dagger 
{\cal A}{\bf x}$. If ${\cal A}$ can be factored to ${\cal A}={\cal B}^\dagger {\cal B}$. Then 
$\|{\cal B}_N {\bf x}\|_2 \longrightarrow \|{\cal B}{\bf x}\|_2$. 
The following is another consequence of Lemma 5.2. 
 
\proclaim  Lemma 5.3. Let ${\cal B}_N$ be a factorization of ${\cal A}_N$, i.e., ${\cal A}_N={\cal 
B}_N^\dagger {\cal B}_N$. Then $\|{\cal B}_N\|$ is bounded independent of $N$. 

\pf By Lemma 5.2, there exists a constant $C$ such that for $N$ large enough, 
$$
\|{\cal B}_N {\bf x}\|_2^2 = {\bf x}^\dagger {\cal A}_N {\bf x}\le {\bf x}^\dagger {\cal A}{\bf x}+C
=\|{\bf x}\|^2_2 \|A(z)\|_\infty + C. 
$$
Hence, $\displaystyle \|{\cal B}_N\|:=\max_{{\bf x}\in\ell_2\atop \|{\bf x}\|_2=1} 
\|{\cal B}_N{\bf x}\|_2$ is bounded. \eop

Note that all ${\cal B}_N$ are  banded upper triangular matrices with one half the band width as that of ${\cal A}$. Thus, 
each row (or column) of ${\cal B}_N$ has finitely many nonzero entries.  Lemma 5.3 
implies that each row (or column) of ${\cal B}_N$ is bounded in $\ell_2$ norm and hence each
entry in any row is bounded. Therefore there exists a  subsequence of ${\cal B}_{N_j}$ such that each
entry with indices $(j,k)$ in ${\cal B}_{N_i}$ converges as $i\longrightarrow +\infty$. 
That is, for any vector ${\bf x}=(x_i)_{i\in {\bf z}}\in 
\ell_2$ with finitely many nonzero entries $x_i$'s, we have
$$
{\cal B}_{N_i}{\bf x} \longrightarrow {\cal B} {\bf x}.
$$
for a bi-infinite matrix ${\cal B}$.  By Lemma 5.2 again, we have ${\bf x}^\dagger 
{\cal B}^\dagger {\cal B}{\bf y}
={\bf x}^\dagger {\cal A}{\bf y}$ for all vectors ${\bf x}$ and ${\bf y}$ with finitely many nonzero entries. However since these are dense in $\ell_2$ we find ${\cal B}^\dagger {\cal B}= {\cal A}$. 
Note that ${\cal B}$ is an upper triangular matrix with the half the band width 
as that of ${\cal A}$. If ${\cal B}$ is a
Toeplitz matrix, we immediately know that $A(z)$ has a factorization such that $A(z)=B(z)^*B(z)$. 
Therefore, we end with

\proclaim Theorem 5.4. Let $P(z,w)$ be a nonnegative Laurent polynomial with no simple zeros. 
Let ${\cal P}$ be 
a bi-infinite Toeplitz matrix with Laurent polynomial entries in $w$, $P_N$ be the central section as described above and ${\cal B_N}$ be its upper triangular Cholesky factor. Then there is a subsequence of  ${\cal B}_N$ convergent to ${\cal B}$ entrywise, where ${\cal P}={\cal B}^\dagger{\cal B}$. 
If ${\cal B}$ is Block Toeplitz, then $P(z,w)$ can be factored into a sum of square 
magnitudes of finitely many polynomials in $z$ and $w$.

Theorem 5.4 provides a computational method to check if a nonnegative Laurent polynomial $P(z,w)$ can
be factorized. That is, we compute Cholesky factorization of central sections of the bi-infinite
Toeplitz matrix ${\cal P}$ associated with $P(z,w)$ and observe if the factorization matrices 
converge to a Toeplitz matrix or not. If they converge, $P(z,w)$ has a factorization.

\sect{6}{Remarks}

\item{1.} It is interesting to point out that the symbol approximation method discussed
in [Lai'94] is very much like the 
Bauer method invented in 1955 (see [Sayed and Kailath'01] and its references).  A slight difference is that  the singular
value decomposition (SVD) instead of the Cholesky decomposition is used to factorize the matrices 
associated with Laurent polynomial $P(z)\ge 0$. 

\item{2.} When $P(z)$ is a matrix polynomial in the univariate setting
[Hardin, Hogan and Sun'04] have demonstrated 
a constructive method to factor $P(z)=Q(z)^\dagger Q(z)$ when $P(z)$ has a nonzero monomial determinant. 

\item{3.} When $P(z)$ is a matrix polynomial in the univariate setting
[Youla and Kazanjian'78] used a Bauer type method to compute the
factorization of $P(z)$. They showed that the solution obtained from the Bauer type method converges to the factorization 
under a weaker condition that 
$$
{1\over 2\pi}\int_{-\pi}^\pi \log \det P(z) d\theta > -\infty
$$
than the positivity condition of $P(z)$. In our Theorem 3.1. we showed the exponential convergence of 
the method which greatly improves their convergence analysis.

\bigskip
\centerline{\bf References}

\ref S. Basu, A constructive algorithm for 2D spectral factorization with rational spectral
factors, IEEE Trans. on Circuits and Systems, 47(2000), 1309--1318.

\ref A. Calderon and R. Pepinsky, On the phases of Fourier coefficients for positive real 
periodic functions, Computing Methods and Phase Problem in X-Ray Crystal Analysis, edited by 
R. Pepinsky, 1952, pp. 339--346.

\ref  I. Daubechies, {\sl Ten Lectures on Wavelets}, SIAM Publications, 
Philadelphia, 1992.

\ref M.  A. Dritschel, On factorization of trigonometric polynomials, Integral Equations and
Operator Theory, 49(2004), 11--42.

\ref J. Geronimo and H. J. Woerdeman, Positive existions, Fej\'er-Riesz factorization and
autoregressive filters in two variables, Annual Math., to appear in 2004. 

\ref D. Hardin, T. Hogen, and Q. Sun, The matrix-valued Riesz Lemma and local orthonormal
bases in shift-invariant spaces, Adv. Comput. Math. 20(2004), 367--384. 

\ref H. Helson, {\sl Lectures on invariant subspaces}, Acedemic Press, New
York, 1964

\ref L. Fej\'er, \"Uber trigonometrische Polynome, J. Reine und Angewandte Mathematik 146(1915), 
53--82.

\ref M. J. Lai, On the computation of Battle-Lemarie's wavelets, Math. Comp. 63(1994), 689--699.

\ref M. J. Lai and J. St\"ockler, Construction of multivariate compactly supported 
tight wavelet frames, submitted,  2004. 

\ref J. W. McLean and H. J. Woerdeman, Spectral factorizations and sums of squares representation
via semi-definite programming, SIAM J. Matrix Anal. Appl. 23(2001), pp. 646--655. 

\ref F. Riesz, \"Uberr ein Problem des Herrn Carath\'eodory, J. Reine und Angewandte Mathematik 
146(1915), 83--87.

\ref W. Rudin, The existence problem for positive definite functions, Illinois J. Math., 
7(1963),pp. 532--539.

\ref A. H. Sayed and T. Kailath, A survey of spetral factorization methods, Numer. Linear Algebra
with Applications 8(2001), 467--496.

\ref D. Youla and N. Kazanjian, Bauer-type factorization of positive matrices and the theory of
matrix polynomials orthogonal on the unit circle, IEEE Trans. Circ. Systems 25(1978), 57--69. 

\end

Thus it is clear that for $x=[x_1, 0, \cdots, 0]^\dagger $,
$x^*\widetilde{P_m}(w)x ={1\over m+1} x_1^* p_0(w) x_1>0$.
For $x=[x_1, x_2, 0, \cdots, 0]^\dagger $, 
$$\eqalign{
x^*\widetilde{P_m}(w)x &= {1\over m+1}(x_1^*p_0x_1+x_2^*p_0x_2)+{1\over m} 
(x_1^* p_1x_2 + x_2^* p_{-1}x_1)\cr
&\ge {1\over m+1}\left( x_1^*p_0x_1+x_2^*p_0x_2 + {m+1\over m}(x_1^* p_1 x_2+ x_2^*p_{-1}x_1\right) \cr
& ={1\over m+1}
[x_1^*, x_2^*]\left[\matrix{p_0(w) & {m+1\over m}p_{-1}(w)\cr {m+1\over m}p_1(w) & p_0(w)\cr}
\right]\left[\matrix{x_1\cr x_2\cr}\right]\cr 
&={1\over m+1} 
[x_1^*, x_2^*]{\cal P}_2 \left[\matrix{x_1\cr x_2\cr}\right]\cr
&\quad + {1\over m+1}
[x_1^*, x_2^*]\left[\matrix{0 & {1\over m}p_{-1}(w)\cr {1\over m}p_1(w) & 0\cr}
\right]\left[\matrix{x_1\cr x_2\cr}\right].\cr }
$$
$P(z,w)\ge \epsilon>0$ implies that 
$$
[x_1^*, x_2^*]{\cal P}_2 \left[\matrix{x_1\cr x_2\cr}\right]\ge \epsilon (|x_1|^2+|x_2|^2). 
$$
Let $C=\max\{\|p_j(w)\|_\infty, j=-n, \cdots, j\}$. Note that 
$$
\eqalign{
&[x_1^*, x_2^*]\left[\matrix{0 & {1\over m}p_{-1}(w)\cr {1\over m}p_1(w) & 0\cr}
\right]\left[\matrix{x_1\cr x_2\cr}\right] \cr
\le &  [x_1, x_2] \left[\matrix{{1\over m}p_{-1}(w)x_2\cr {1\over m}p_1(w)x_1}\right] \cr
\le & {C\over m} (2x_1x_2) \le {C\over m}(|x_1|^2+|x_2|).}
$$
If we choose $m$ so large that ${C\over m} <\epsilon$, then we will have
$x^*\widetilde{P_m}(w)x >0$. 

For $x=[x_1,x_2,x_3, 0, \cdots, 0]^\dagger $, 
$$\eqalign{
&x^*\widetilde{P_m}(w)x \cr 
&= {1\over m+1} [x_1^*,x_2^*,x_3^*] 
\left[\matrix{ p_0(w) & {m+1\over m}p_{-1}(w) & {m+1\over m-1}p_{-2}\cr 
{m+1\over m}p_1(w) & p_0(w) & {m+1\over m}p_{-1}(w)\cr
{m+1\over m-1}p_2(w) & {m+1\over m}p_1(w) & p_0(w)\cr}\right] \left[\matrix{x_1\cr x_2\cr x_3\cr}
\right]\cr
&= {1\over m+1} [x_1^*,x_2^*,x_3^*] {\cal P}_3 \left[\matrix{x_1\cr x_2\cr x_3\cr}
\right]\cr
& \qquad + {1\over m+1} [x_1^*,x_2^*,x_3^*] 
\left[\matrix{ 0 & {1\over m}p_{-1}(w) & {2\over m-1}p_{-2}\cr 
{1\over m}p_1(w) & 0 & {1\over m}p_{-1}(w)\cr
{2\over m-1}p_2(w) & {1\over m}p_1(w) & 0\cr}\right] \left[\matrix{x_1\cr x_2\cr x_3\cr}
\right]\cr
&\ge {\epsilon\over m+1} \|[x_1,x_2,x_3]\|^2  \cr
&\quad  + {1\over m+1} [x_1^*,x_2^*,x_3^*] 
\left[\matrix{ 0 & {1\over m}p_{-1}(w) & {2\over m-1}p_{-2}\cr 
{1\over m}p_1(w) & 0 & {1\over m}p_{-1}(w)\cr
{2\over m-1}p_2(w) & {1\over m}p_1(w) & 0\cr}\right] \left[\matrix{x_1\cr x_2\cr x_3\cr}
\right].\cr}
$$
Note that we have
$$
\eqalign{
& \|\left[\matrix{ 0 & {1\over m}p_{-1}(w) & {2\over m-1}p_{-2}\cr 
{1\over m}p_1(w) & 0 & {1\over m}p_{-1}(w)\cr
{2\over m-1}p_2(w) & {1\over m}p_1(w) & 0\cr}\right] 
 \left[\matrix{x_1\cr x_2\cr x_3\cr}
\right] \|_2^2 \cr
&\le \left({2C\over m-1}\right)^2(x_2^2+x_3^2+x_1^2+x_3^2+x_1^2+x_2^2) \cr 
&\le \sqrt{2}\left({2C\over m-1}\right)^2(|x_1|^2+|x_2|^2+|x_3|^2).\cr}
$$
It follows that 
$$
\|\left[\matrix{ 0 & {1\over m}p_{-1}(w) & {2\over m-1}p_{-2}\cr 
{1\over m}p_1(w) & 0 & {1\over m}p_{-1}(w)\cr
{2\over m-1}p_2(w) & {1\over m}p_1(w) & 0\cr}\right] \|_2 \le {4C\over m-1}
$$
and 
$$
x^*\widetilde{P_m}(w)x\ge {\epsilon\over m+1} \|[x_1,x_2,x_3]\|^2  -{1\over m+1}{2\sqrt{2}C\over m-1}.
$$
If we choose $m$ so large to make $2\sqrt{2}C/(m-1) <\epsilon$, we have
$x^* \widetilde{P_m}(w)x>0$ in this case.